\documentclass{article}                                                                           
\usepackage[margin=1.3in]{geometry}                                                     
\usepackage{setspace}                                                                             
\usepackage{indentfirst}      
\usepackage[table]{xcolor}                                                                    
\renewcommand\thesection{\Roman{section}}                                         
\renewcommand\thesubsection{\thesection.\Alph{subsection}}                   
\renewcommand\thesubsubsection{\thesubsection.\arabic{subsubsection}} 
\usepackage[explicit]{titlesec}                                                                  
\titleformat{\section}{\bfseries}{\thesection.}{1em}{\MakeUppercase{#1}}  
\titleformat{\subsection}{\bfseries}{\thesubsection.}{1em}{#1}                   
\titleformat{\subsubsection}{\itshape}{\thesubsubsection.}{1em}{#1}          
\usepackage{lastpage} 
\usepackage{color,soul}                                                                            
\usepackage{amsmath,array}                                                                            
\usepackage[]{footmisc}                                                                          
                                  
\usepackage{graphicx}                                                                             
\usepackage{multirow}                                                                             
\usepackage{caption}                                                                              
\usepackage[labelformat=simple]{subcaption}                                           
\usepackage[colorlinks=true, allcolors=blue]{hyperref}

\title{Conditions for Translation and Scaling Invariance of the Neutron Diffusion Equation}
\author{
\\Jesse F.~Giron$^{\text{a},\text{b}}$\footnote{Email: \href{mailto:jgiron@lanl.gov}{jgiron@lanl.gov}} , Scott  D.~Ramsey$^{\text{a}}$\footnote{Email: \href{mailto:ramsey@lanl.gov}{ramsey@lanl.gov}} , and Brian A.~Temple$^{\text{a}}$\footnote{Email: \href{mailto:temple@lanl.gov}{temple@lanl.gov}} \\[4pt]
LA-UR 18-20828\\[4pt] 
\textit{$^a$Los Alamos National Laboratory, Applied Physics}\\[-10pt]       
\textit{P.O. Box 1663, MS P225, Los Alamos, New Mexico 87545} \\[-2pt]
\textit{$^b$Department of Physics, Box 871504}\\[-10pt]       
\textit{Arizona State University, Tempe, Arizona  85287-1504} \\[-2pt]}
\date{}
\usepackage{graphicx} 
\usepackage{booktabs} 
\usepackage{microtype} 

\begin{document}
\maketitle
\pagebreak
\begin{abstract}
Lie group methods are applied to the time-dependent, monoenergetic neutron diffusion equation in materials with spatial and time dependence. To accomplish this objective, the underlying 2nd order partial differential equation (PDE) is recast as an exterior differential system so as to leverage the isovector symmetry analysis approach. Some of the advantages of this method as compared to traditional symmetry analysis approaches are revealed through its use in the context of a 2nd order PDE. In this context, various material properties appearing in the mathematical model (e.g., a diffusion coefficient and macroscopic cross section data) are left as arbitrary functions of space and time. The symmetry analysis that follows is restricted to a search for translation and scaling symmetries; consequently the Lie derivative yields specific material conditions that must be satisfied in order to maintain the presence of these important similarity transformations.
The principal outcome of this work is thus the determination of analytic material property functions that enable the presence of various translation and scaling symmetries within the time- dependent, monoenergetic neutron diffusion equation. The results of this exercise encapsulate and generalize many existing results already appearing in the literature. While the results contained in this work are primarily useful as phenomenological guides pertaining to the symmetry behavior of the neutron diffusion equation under certain assumptions, they may eventually be useful in the construction of exact solutions to the underlying mathematical model. The results of this work are also useful as a starting point or framework for future symmetry analysis studies pertaining to the neutron transport equation and its many surrogates.
\end{abstract}
\pagebreak
\section{INTRODUCTION}
Neutron diffusion equations are pervasive in the nuclear engineering community, from their role in introducing students to the phenomenology of neutral particle transport \cite{lamarsh}, to their widespread use in the reactor design community. In the motivation, development, and application of these equations a variety of approximate solution techniques are employed, including multi-group energy discretizations and a wide variety of space and time differencing techniques~\cite{clark_hansen}\cite{henry}. In many cases, systems of neutron diffusion equations are discretized and solved as large systems of algebraic equations using computer codes.

Given these developments, the role and usefulness of closed-form or semi-analytic solutions to neutron diffusion equations would seem to be relegated to the classroom as practice tools for aspiring nuclear engineers. However, as demonstrated in a variety of other contexts including solid mechanics, fluid flow, heat transport, and wave propagation \cite{guenther_lee}, exact solutions of partial differential equations can still serve an important purpose even in fields increasingly dominated by computational horsepower. In particular, code verification via ``exact solutions as test problems" is becoming an increasingly visible and important aspect of the code development paradigm~\cite{oberkampf_roy}\cite{roache}. To this end, several compendia of such test problems fit for the nuclear engineering context have been developed by many authors~\cite{ganapol}-\cite{trahan_sweezy_giron}.

While the current state of solution methods for the governing partial differential equations (PDEs) may otherwise seem to consist of various ad hoc methods~\cite{guenther_lee}, the unified theory for achieving exact solutions or reduced-order structures (which will be more amenable to high-accuracy numerical solution techniques) is found in the field of symmetry analysis, also variously known as group-theoretic or Lie group techniques \cite{bluman_anco}-\cite{stephani} . In short, these methods provide a unified setting through which to determine the symmetries inherent to PDEs (and, more generally, other algebraic, differential, integral, or discrete structures); if these can be shown to exist, they enable changes of variables through which an order reduction, or sometimes an exact solution can be obtained.

Moreover, symmetry analysis techniques as applied to differential equations are also valuable from the standpoint of developing or reinforcing physical intuition. Some of the simplest but most useful symmetries identifiable through Lie group techniques include physical transformations or similarities such as translations, scalings, rotations, and projections. If a differential equation motivated by a physical process is found to have one or more of these symmetries, broad statements regarding its applicability can be constructed. For example, differential equations possessing scaling symmetries can often be shown to be valid across temporal, spatial, or other scales; similarly, translation or rotation-invariant equations can be expected to be valid across various orientations. Developing an understanding of these phenomena is important not only from the standpoint of constructing exact solutions with special properties, but also in terms of guiding or interpreting experimental or computational activities. As developed extensively by Barenblatt~\cite{barenblatt}, physically intuitive symmetries often give rise to the quasi-limiting phenomenon known as ``intermediate asymptotics", where much of the essential physics inherent to a problem of interest is captured by similarity processes.

Following the initial revitalization of Lie's original techniques in the 1950s by Birkhoff~\cite{birkhoff} and Ovsiannkov~\cite{ovsiannkov}, the following years have been witness to a veritable explosion in the use of symmetry methods in a wide variety of physical models. For example, in the context of neutron transport and its surrogate models, symmetry analysis has been performed on:
\begin{itemize}
\item A time-independent formulation of neutron diffusion, by Tsyfra and Czyzycki~\cite{tsyfra_czyzycki},

\item Group-invariant differencing schemes by Axford~\cite{axford1}-\cite{axford_noh}, Grove~\cite{grove}, Jaegers~\cite{jaegers}, and Melenshko~\cite{melenshko}.

\item Various instantiations of separable time-dependent neutron diffusion equations by Ramsey, et al.~\cite{ramsey_tellez}.
\end{itemize}
Similarly, exhaustive studies of the time-dependent heat conduction equation have been performed by numerous authors~\cite{bluman_anco}-\cite{bluman_kumei}\cite{olver} and are useful as guiding studies due to the correspondence between heat conduction and neutral particle diffusion.

The current study intends to initiate a broad application of symmetry analysis techniques in the context of the neutron transport equation and its surrogates. To draw a connection to previous work, emphasis will first be placed on a diffusion model. To this end, in this work a monoenergetic diffusion equation with space and time-dependent material properties will be developed and analyzed for the existence of translational and scaling symmetry groups (these being two common sub-classes of similarity transformations). The purpose of these calculations is to determine a correspondence between the functional forms material properties may assume and the resulting translation and scaling transformation properties they allow. These outcomes are intended to further a geometric and physical interpretion of the governing equation as discussed above; the construction of exact solutions resulting from these symmetries will be relegated to a future work. 

In support of these developments, Sec. II includes an overview of the attendant mathematical model, including a brief review of symmetry analysis methods. Sec. III contains the calculation of the admissible translation and scaling groups for the diffusion equation under consideration, and their connection to the admissible functional forms of the included material properties. A discussion of these results and their connections to other bodies of work is provided in Sec. IV. Finally conclusions and recommendations for future study follow in Sec. V.
\section{MATHEMATICAL MODEL}\label{sec:two}
\subsection{Neutron Diffusion}

The traditional neutron diffusion equation may be viewed as a representation of the neutron transport equation in a low-order spherical harmonics approximation. Derivations of this equation are provided by numerous authors, including Lamarsh \cite{lamarsh}, Duderstadt and Hamilton \cite{duderstadt_hamilton}, and many others; it may be written in generalized 1D curvilinear coordinates as

\begin{equation}\label{timedependent_diff_eq}
\frac{1}{\rm{v}}\frac{\partial \phi}{\partial t} = \frac{1}{r^n}\frac{\partial}{\partial r}\bigg[r^n D(r,t)\frac{\partial \phi}{\partial r}\bigg]+\bigg[\bar{\nu}\Sigma_f(r,t)-\Sigma_a(r,t)\bigg]\phi,
\end{equation}
where $\phi$ is the scalar neutron flux interpreted as a function of time $t$ and the general 1D spatial coordinate $r$: $n = 0,1,2$ for Cartesian, cylindrical, and spherical coordinates, respectively, and $\rm{v}$ is the neutron speed. The medium through which neutrons are diffusing is characterized by a diffusion coefficient $D$ and a macroscopic total absorption cross section $\Sigma_a$; as fission processes are included, so must be a macroscopic fission cross section $\Sigma_f$ and average number of neutrons released per fission $\bar{\nu}$. These material properties are regarded as unspecified functions of time and space, as indicated in Eq. (\ref{timedependent_diff_eq}).

As written, Eq. (\ref{timedependent_diff_eq}) represents a ``monoenergetic" model. The material properties appearing in Eq. (\ref{timedependent_diff_eq}) are implicitly taken as averages with respect to a prescribed neutron energy spectrum, and the energy dependence of $\phi$ is not explicitly represented; the scalar flux is assumed to be evaluated at either a discrete or averaged value of the neutron energy. While Eq. (\ref{timedependent_diff_eq}) may be of limited practical utility, the forthcoming analysis conducted in its context will reveal many concepts that may eventually be extended to more complicated and realistic models (e.g., multi-group $P_N$ or $S_N$ approximations). 

\subsection{Differential Forms}
\label{sec:forms}

As noted in Sec. II. A., the goal of this work is to subject Eq.~(\ref{timedependent_diff_eq}) to a rigorous program of symmetry analysis, in the interest of uncovering its invariance properties with respect to scaling and translation transformations. The systematic means through which to conduct symmetry analysis is discussed by many authors including Olver~\cite{olver}, Bluman and collaborators~\cite{bluman_anco}-\cite{bluman_kumei}, Cantwell~\cite{cantwell}, Stephani~\cite{stephani}, and many others. Inherent to this procedure is a geometric interpretation of the equation(s) under investigation. In the classical symmetry analysis methods outlined by the aforementioned authors, this interpretation manifests itself through the identification of the relevant differential equations as purely algebraic structures in a suitable higher-dimensional manifold. After an invariance analysis is conducted, relationships between independent variables, dependent variables, and their derivatives are recovered via prolongation formulae that arise from requiring invariance of the definition of derivatives across the coordinate transformations under investigation.

In 1971, Harrison and Estabrook~\cite{harrison_estabrook} formulated an alternative approach to classical symmetry analysis methods using the language of differential forms. While Harrison and Estabrook~\cite{harrison_estabrook} and many subsequent authors~\cite{edelen}\cite{suhubi} were able to demonstrate that the symmetry analysis results arising from their approach were identical to those found using the classical methods, the differential form or \emph{isovector} method provides an elegant, intuitive geometric setting from which one can initiate a symmetry analysis of differential equations. From a practical standpoint, given the well-established properties of exterior and Lie derivatives (see Sec. II.C.), the need for sometimes cumbersome prolongation formulae (especially in the context of higher-order differential equations) appears to be obviated.

Perhaps the only ``drawback" to the isovector approach is the necessity of casting the relevant differential equations as a first order system (a requirement for the further re-casting of that system as an exterior differential system, or EDS). While this procedure may prove tedious for high-order systems (or impossible for certain others), the complication is nearly trivial for a single second-order ODE. In the context of Eq.~(\ref{timedependent_diff_eq}), the equivalent ``first-order" system follows from the substitution

\begin{equation}\label{eqn:w_def}
w\left(r,t\right) = \frac{\partial \phi}{\partial r},
\end{equation}
so that Eq. (\ref{timedependent_diff_eq}) becomes
\begin{equation}\label{eqn:diffeq_final}
\frac{1}{\rm{v}}\frac{\partial \phi}{\partial t} = \frac{1}{r^n}\frac{\partial}{\partial r}\bigg[r^n D(r,t)w\bigg]+\bigg[\bar{\nu}\Sigma_f(r,t)-\Sigma_a(r,t)\bigg]\phi.
\end{equation}
Equations (\ref{eqn:w_def}) and (\ref{eqn:diffeq_final}) may then be multiplied by the differential volume element $dt\wedge dr$ to yield 
\begin{eqnarray}\label{eqn:mu1}
\mu_1 = \frac{1}{\rm{v}} d\phi \wedge dr + nr^{-1} D d\phi \wedge dt + D_r d\phi \wedge dt +D dw\wedge dt
	+ \big(\bar{\nu}\Sigma_f - \Sigma_a\big)\phi dt\wedge dr,
\end{eqnarray}
and 
\begin{eqnarray}\label{eqn:mu_2}
\mu_2 = w dt\wedge dr + d\phi\wedge dt,
\end{eqnarray}
where Eqs. (\ref{eqn:mu1}) and (\ref{eqn:mu_2}) are referred to as a system of 2-forms. Note, as we proceed we define $\Gamma(r,t) \equiv \bar{\nu}\Sigma_f(r,t) - \Sigma_a(r,t)$, which will be known as the ``gamma coefficient". In addition, note $D_r = \frac{\partial D}{\partial r}$. In Eqs. (\ref{eqn:mu1}) and (\ref{eqn:mu_2}), the operator $d$ is referred to as an exterior derivative, and the operator $\wedge$ used to multiply differentials is known as a wedge product. The salient properties are
\begin{eqnarray}\label{antisymmetric_prop}
dq \wedge dp = - dp \wedge dq,
\end{eqnarray}
and 
\begin{eqnarray}\label{null_prop}
dq\wedge dq = 0
\end{eqnarray}
for all general coordinates $p$ and $q$ (i.e. wedge-multiplication of differentials is antisymmetric). More comprehensive overviews of differential geometry and the structures appearing therein are provided by Edelen~\cite{edelen}, Suhubi~\cite{suhubi}, Bourbaki~\cite{bourbaki}, and Albright et. al.~\cite{jason} .

As written, Eqs.~(\ref{eqn:mu1}) and (\ref{eqn:mu_2}) indicate that $r$, $t$, $\phi$, and $w$ are entirely indpendent of each other, and thus $\mu_1$ and $\mu_2$ represent differential objects in an appropriately expanded space. To establish the equivalence between these structures and Eqs.~(\ref{eqn:w_def}) and (\ref{eqn:diffeq_final}), independent and dependent variable relationships are once again enforced (a process referred to as ``sectioning" by Harrison and Estabrook). Letting $\phi = \phi\left(r,t\right)$ and $w = w\left(r,t\right)$, the exterior derivatives (or total differentials) appearing in Eqs.~(\ref{eqn:mu1}) and (\ref{eqn:mu_2}) can be expanded to yield

\begin{eqnarray}\label{flux_differential_expansion}
\frac{1}{\rm{v}} \bigg(\frac{\partial \phi}{\partial r}dr + \frac{\partial \phi}{\partial t}dt\bigg)\wedge dr + nr^{-1} D \bigg(\frac{\partial \phi}{\partial r}dr + \frac{\partial \phi}{\partial t}dt\bigg)\wedge dt 
	+ D_r\bigg(\frac{\partial \phi}{\partial r}dr + \frac{\partial \phi}{\partial t}dt\bigg)\wedge dt \nonumber\\+D\bigg( \frac{\partial w}{\partial r}dr + \frac{\partial w}{\partial t}dt\bigg)\wedge dt
	+ \Gamma\phi\; dt\wedge dr
\end{eqnarray}
and
\begin{eqnarray}\label{w_differential_expansion}
w dt\wedge dr + \bigg(\frac{\partial \phi}{\partial r}dr + \frac{\partial \phi}{\partial t}dt\bigg)\wedge dt .
\end{eqnarray}
Using the Eqs. (\ref{antisymmetric_prop}) and (\ref{null_prop}), we obtain the following:
\begin{eqnarray}\label{use_antisymmetry_flux}
\bigg[\frac{-1}{\rm{v}} \frac{\partial \phi}{\partial t} + nr^{-1} D\frac{\partial \phi}{\partial r}+ D_r\frac{\partial \phi}{\partial r} +D \frac{\partial w}{\partial r} + \Gamma\phi\bigg]\;dt\wedge dr,
\end{eqnarray}
and
\begin{eqnarray}\label{use_antisymmetry_w}
\bigg(w - \frac{\partial \phi}{\partial r}\bigg)\;dt\wedge dr .
\end{eqnarray}
Setting these relations equal to zero (a process referred to as ``annulling" by Harrison and Estabrook~\cite{harrison_estabrook}), the nontrivial solution that follows is given by Eqs.~(\ref{eqn:w_def}) and (\ref{eqn:diffeq_final}).

Equations (\ref{eqn:mu1}) and (\ref{eqn:mu_2}) will be used for the symmetry analysis studies to follow, given their underlying equivalence to Eq. (\ref{timedependent_diff_eq}), the mathematical model of principal interest to this study.

\subsection{Invariance}

Symmetry analysis techniques as applied to the study of differential equations are in a sense a unification theory for disparate, otherwise ad hoc methods for the solution of those equations. The key point surrounding these techniques is that if a differential equation possesses symmetries, they will enable a change of coordinates through which the original equation may either be reduced to a simpler structure (e.g., from a PDE to an ODE) or solved outright. Of interest to this work are the scenarios wherein Eq. (\ref{timedependent_diff_eq}) - or more precisely, its EDS representation given by Eqs.~(\ref{eqn:mu1}) and (\ref{eqn:mu_2}) - is invariant under continuous scaling and translation transformations. The resulting reduced-order structures and/or solutions will then be characterized by translation and scale invariance; in the latter case, this extends their validity across dimensional scale (i.e., the solution is expected to be valid in any set of units). This phenomenon enables the use of scaled experiments, etc., to diagnose situations that may otherwise be prohibitive from a cost or other standpoint.

To proceed, the objective is to determine for what values of the constants $\tau_1$-$\tau_4$ and $s_1$-$s_4$ the global transformations given by
\begin{eqnarray}
r_{\rm{new}} &=& \epsilon a_1 + e^{\epsilon a_2} r ,	\label{eqn:rnew} \\
t_{\rm{new}} &=& \epsilon a_3 + e^{\epsilon a_4} t ,	\label{eqn:tnew} \\
\phi_{\rm{new}} &=& \epsilon a_5 + e^{\epsilon a_6} \phi ,	\label{eqn:phinew} \\
w_{\rm{new}} &=& \epsilon a_7 + e^{\epsilon a_8} w ,	\label{eqn:wnew}
\end{eqnarray}
leave Eqs.~(\ref{eqn:mu1}) and (\ref{eqn:mu_2}) unchanged in structure; that is, when the substitutions given by Eqs.~(\ref{eqn:rnew})-(\ref{eqn:wnew}) are put into Eqs.~(\ref{eqn:mu1}) and (\ref{eqn:mu_2}), the resulting relations are unchanged aside from the indexing $r \to r_{\rm{new}}$, $t \to t_{\rm{new}}$, and so forth. For example, in the context of Sec. II.B.,
\begin{eqnarray}
\mu_1\left(r_{\rm{new}},t_{\rm{new}},\phi_{\rm{new}},w_{\rm{new}}\right) &=& \mu_1\left(r,t,\phi,w\right), \label{eqn:mu1_invariance} \\
\mu_2\left(r_{\rm{new}},t_{\rm{new}},\phi_{\rm{new}},w_{\rm{new}}\right) &=& \mu_2\left(r,t,\phi,w\right). \label{eqn:mu2_invariance}
\end{eqnarray}

This global concept of invariance may be alternatively realized on the local level in terms of an appropriate vector field. To see this, a Taylor-series expansion of the left-hand side of Eqs.~(\ref{eqn:mu1_invariance}) and (\ref{eqn:mu2_invariance}) about the identify element $\epsilon = 0$ of the transformations given by Eqs.~(\ref{eqn:rnew})-(\ref{eqn:wnew}) yields
\begin{equation}\label{taylor_expansion_mu_new}
\mu_{i,\rm{new}} =\mu_i +\left.\epsilon\frac{\partial\mu_i}{\partial\epsilon}\right|_{\epsilon=0} + \left. \frac{1}{2}\epsilon^2 \frac{\partial^2\mu_i}{\partial\epsilon^2}\right|_{\epsilon=0} + \cdots
\end{equation}
The various derivatives appearing in the Taylor-series expansion may be expanded using the chain rule
\begin{equation}
\frac{\partial}{\partial\epsilon} = \left.\frac{\partial r_{\rm{new}}}{\partial\epsilon}\right|_{\epsilon=0}\frac{\partial}{\partial r} + \left.\frac{\partial t_{\rm{new}}}{\partial\epsilon}\right|_{\epsilon=0}\frac{\partial}{\partial t}+\left.\frac{\partial \phi_{\rm{new}}}{\partial\epsilon}\right|_{\epsilon=0}\frac{\partial}{\partial \phi}+\left.\frac{\partial w_{\rm{new}}}{\partial\epsilon}\right|_{\epsilon=0}\frac{\partial}{\partial w},
\end{equation}
or, with Eq.~(\ref{eqn:rnew})-(\ref{eqn:wnew}), 
\begin{equation}\label{epsilon_group_generator}
\frac{\partial}{\partial\epsilon}=\chi = (a_1 + a_2r)\frac{\partial}{\partial r}+(a_3 + a_4t)\frac{\partial}{\partial t} + (a_5 + a_6\phi)\frac{\partial}{\partial \phi}
	+ (a_7 + a_8w)\frac{\partial}{\partial w}.
\end{equation}
With Eq.~(\ref{epsilon_group_generator}), Eq.~(\ref{taylor_expansion_mu_new}) becomes
\begin{equation}
\mu_{i,\rm{new}} = \mu_i+\epsilon\chi\mu_i+\frac{1}{2}\epsilon^2 \chi\chi\mu_i +\cdots
\end{equation}
and the invariance condition given by Eq.~(\ref{eqn:mu1_invariance}) and (\ref{eqn:mu2_invariance})
becomes
\begin{equation}\label{taylor_invariance}
\epsilon\chi\mu_i + \frac{1}{2}\epsilon^2\chi\chi\mu_i +\cdots =0.
\end{equation}
The nontrivial solution of Eq.~(\ref{taylor_invariance}) is $\chi\mu_i =0$, or 
\begin{eqnarray}\label{group_generator}
(a_1 + a_2r)\frac{\partial \mu_i}{\partial r}+(a_3 + a_4t)\frac{\partial \mu_i}{\partial t} + (a_5 + a_6\phi)\frac{\partial \mu_i}{\partial \phi}
	+ (a_7 + a_8w)\frac{\partial \mu_i}{\partial w} = 0,
\end{eqnarray}
when
\begin{equation}\label{eqn:mu_i_zero}
\mu_i = 0.
\end{equation}
This localized or infinitesimal invariance condition is entirely equivalent to the global invariance criterion provided in Eqs.~(\ref{eqn:mu1_invariance}) and (\ref{eqn:mu2_invariance}).

The operator $\chi$ appearing in Eq. (\ref{group_generator}) is variously referred to as the vector field generated by the group of transformations given in Eqs.~(\ref{eqn:rnew})-(\ref{eqn:wnew}), or more simply, the ``group generator". As noted by Olver~\cite{olver} and Stephani~\cite{stephani}, it is also an example of a Lie derivative, or a generalization of the more familiar directional derivative operator as appearing in vector calculus. In this case, the Lie derivative is defined on the higher-dimensional manifold spanned by all independent and dependent variables associated with a system under investigation. The interpretation of $\chi$ as a Lie derivative also makes clear the advantages of a differential geometric interpetation of symmetry analysis of differential equations: in particular, in the evaluation of Eq.~(\ref{eqn:mu1}), structures such as
\begin{equation}
\chi d \left( q \right)
\end{equation}
frequently appear. Their evaluation is facilitated by the commutability of the exterior and Lie derivative operations (see, for example, Edelen~\cite{edelen}):
\begin{equation}\label{commutability}
\chi (d q) = d (\chi q )
\end{equation}
for any function $q$. Equation (\ref{commutability}) will be of central importance in the evaluation of Eq. (\ref{eqn:mu1}). 
\section{SYMMETRY ANALYSIS OF THE NEUTRON DIFFUSION EQUATION}
\label{sec:calculations}

Equations (\ref{eqn:mu1}) and (\ref{eqn:mu_2}) are invariant under combined translation and scaling transformations, provided Eqs.~(\ref{group_generator}) and (\ref{eqn:mu_i_zero}) are satisfied. Expanding Eq. (\ref{group_generator}), using the property from Eq. (\ref{commutability}), with $i=1$ yields:
\begin{multline}
\frac{1}{\rm{v}}\left\{\chi(r)d\phi\wedge dr + r d\left[\chi(\phi)\right]\wedge dr + r d\phi\wedge d\left[\chi(r)\right]\right\}+n\left\{\chi(D)d\phi\wedge dt + D d\left[\chi\left(\phi\right)\right]\wedge dt + Dd\phi\wedge d\left[\chi\left(t\right)\right] \right\}\\
+\left\{\chi\left(r\right)D_r d\phi\wedge dt +r\chi\left(D_r\right)d\phi\wedge dt +rD_r d\left[\chi\left(\phi\right) \right]\wedge dt + rD_r d\phi\wedge d\left[\chi\left(t\right)\right] \right\}\\
+\left\{ \chi(r)D dw\wedge dt + r\chi(D) dw\wedge dt + rD d\left[\chi(w)\right]\wedge dt + rD dw\wedge d\left[\chi(t)\right] \right\}\\
+\left\{ \chi(r)\Gamma\phi dt\wedge dr + r \chi(\Gamma)\phi dt\wedge dr + r\Gamma\chi(\phi) dt\wedge dr + r\Gamma\phi d\left[\chi(t)\right]\wedge dr + r\Gamma\phi dt\wedge d\left[\chi(r)\right] \right\}\\
=0.
\end{multline}
Applying each derivative of the group generator leads to
\begin{multline}
\frac{1}{\rm{v}}\left\{\left(a_1 +a_2r\right)d\phi\wedge dr + r a_6 d\phi\wedge dr + r a_2 d\phi \wedge dr \right\} + n \{ \left(a_1 + a_2r\right)D_r d\phi\wedge dt + \left( a_3 + a_4t\right)D_t d\phi\wedge dt\\
+ D a_6 d\phi\wedge dt + D a_4 d\phi\wedge dt \} + \{\left(a_1 + a_2r\right)D_r d\phi\wedge dt + r\left(a_1+a_2r\right)D_{rr} d\phi\wedge dt + r\left(a_3+a_4t\right)D_{rt} d\phi\wedge dt\\
+rD_ra_6 d\phi\wedge dt + rD_ra_4 d\phi\wedge dt\} + \{\left(a_1+a_2r\right)D dw\wedge dt + r\left(a_1+a_2r\right)D_r dw\wedge dt + r\left(a_3+a_4t\right)D_t dw\wedge dt \\
+ rD a_8 dw\wedge dt+ rD a_4 dw\wedge dt\}+ \{\left(a_1+a_2r\right)\Gamma\phi dt\wedge dr + \phi r \left(a_1+a_2r\right)\Gamma_r dt\wedge dr + \phi r \left(a_3+a_4t\right)\Gamma_t dt\wedge dr \\
+r\Gamma\left(a_5+a_6\phi\right) dt \wedge dr + r\Gamma\phi a_4 dt\wedge dr + r\Gamma\phi a_2 dt\wedge dr\}
=0,
\end{multline}
which simplifies to
\begin{multline}\label{balance_form}
\frac{1}{\rm{v}}\left[(a_1 + a_2r)+r(a_6+a_2)\right]d\phi\wedge dr + nD\left(a_6+a_4\right)d\phi\wedge dt
+ n\left[(a_1+a_2r)D_r + (a_3+a_4t)D_t\right]d\phi\wedge dt \\+ D_r\left(a_1+a_2r\right)d\phi\wedge dt
+ r\left[(a_1+a_2r)D_{rr} + (a_3+a_4t)D_{rt}\right]d\phi\wedge dt
+ D_r\left[(a_1+a_2r)+r(a_6+a_4)\right]d\phi\wedge dt\\
+ r\left[(a_1+a_2r)D_r + (a_3+a_4t)D_t\right]dw \wedge dt 
+ D\left[(a_1+a_2r)+r(a_8+a_4)\right]dw\wedge dt\\
+ r\phi\left[(a_1+a_2r)\Gamma_r + (a_3+a_4t)\Gamma_t\right]dt \wedge dr
+ \Gamma\phi\left[(a_1+a_2r)+r(a_4+a_2) \right]dt\wedge dr +\Gamma r\left(a_5+a_6\phi\right)dt\wedge dr\\
=0.
\end{multline}
Note, for convenience, Eq. (\ref{eqn:mu1}) was multiplied by $r$ in order to prevent an undefined function when $r=0$.

Similarly, expanding Eq. (\ref{group_generator}) with $i=2$ yields:
\begin{equation}
\chi(w)dt\wedge dr + w d\left[\chi(t)\right]\wedge dr + w dt\wedge d\left[\chi(r)\right] + d\left[\chi(\phi)\right]\wedge dt + d\phi\wedge d\left[\chi(t)\right] =0.
\end{equation}
Again, applying the appropriate derivatives leads to
\begin{equation}
\left(a_7+a_8w\right)dt\wedge dr + a_4w dt\wedge dr + a_2w dt\wedge dr + a_6 d\phi\wedge dt + a_4d\phi\wedge dt =0,
\end{equation}
which simplifies to
\begin{equation}\label{contact_form}
\left(a_7+a_8w\right)dt\wedge dr + \left(a_4 + a_2\right)w dt \wedge dr + \left(a_6+a_4\right)dt\wedge dr=0.
\end{equation}

The conditions $\mu_1 = 0$ and $\mu_2 = 0$ have been used to eliminate one of the basis 2-forms appearing in each of the expanded relations. These conditions are implicitly included in Eqs.~(\ref{balance_form}) and (\ref{contact_form}).

For the preceding equations to be nontrivially satisfied, they must be identities in each unique basis 2-form appearing within them. From Eq. (\ref{balance_form}):
\begin{itemize}
\item for $d\phi\wedge dr$
\begin{multline}\label{phi_r_contact}
\frac{1}{\rm{v}}\bigg(a_1 +a_2r\bigg)+\frac{r}{\rm{v}}\bigg(a_6+a_2\bigg) - \frac{1}{D \rm{v}}\bigg\{\left(a_1+a_2r\right)D+rD\left(a_8+a_4\right)\\
+r\left[\left(a_1+a_2r\right)D_r + \left(a_3+a_4t\right)D_t\right] \bigg\}=0,
\end{multline}
\item for $d\phi\wedge dt$
\begin{multline}\label{phi_t_contact}
n\bigg\{\big[\left(a_1+a_2r\right)D_r + \left(a_3+a_4t\right)D_t\big]+D\left(a_6+a_4\right)\bigg\}
+D_r\big[\left(a_1+a_2r\right)+a_6r +a_4r\big]\\
+\left(a_1+a_2r\right)D_{rr}+\left(a_3+a_4t\right)D_{rt}
-\left(nD+rD_r\right)\frac{1}{rD}\bigg\{ \left(a_1+a_2r\right)D +r\big[\left(a_1+a_2r\right)D_r\\
+ \left(a_3+a_4t\right)D_t +D\left(a_8+a_4\right)\big]\bigg\}=0,
\end{multline}
\item for $dt \wedge dr$
\begin{multline}\label{r_t_contact}
\left(a_1+a_2r\right)\Gamma \phi +r\phi\big[\left(a_1+a_2r\right)\Gamma_r +\left(a_3+a_4t\right)\Gamma_t\big]
-\frac{\Gamma\phi}{D}\bigg\{\left(a_1+a_2r\right)D + r\big[\left(a_1+a_2r\right)D_r\\ 
+\left(a_3+a_4t\right)D_t+ D\left(a_8+a_4\right)\big]\bigg\}+r\Gamma\bigg[\phi\left(a_2+a_4+a_6\right)+a_5\bigg]=0,
\end{multline}
\end{itemize}
and from Eq. (\ref{balance_form})
\begin{itemize}
\item for $dt\wedge dr$
\begin{equation}\label{t_r_balance}
a_7 = 0,
\end{equation}
\item for $d\phi \wedge dt$
\begin{equation}\label{phi_t_balance}
a_6 -a_2 -a_8 =0.
\end{equation}
\end{itemize}
From these conditions, we find the following determining equations for the material properties:
\begin{equation}\label{D-eq}
(a_1+a_2r)D_r + (a_3+a_4t)D_t +(a_8+a_4-a_2-a_6) D = 0,
\end{equation}
\begin{equation}\label{nD_Dequation}
a_1 D_r = 0,
\end{equation}
\begin{equation}\label{Dr-eq}
(a_1+a_2r)D_{rr} + (a_3+a_4t)D_{rt} +(a_4-a_2) D_r = 0,
\end{equation}
\begin{equation}\label{nDa1equation}
nDa_1 = 0,
\end{equation}
\begin{equation}\label{gamma_eqn}
(a_1+a_2r)\Gamma_r + (a_3+a_4t)\Gamma_t + a_4\Gamma =0,
\end{equation}
\begin{equation}\label{a5-eqn}
a_5 =0.
\end{equation}
In addition to the above, Eqs.~(\ref{t_r_balance}) and (\ref{phi_t_balance}) are also determining equations. Therefore, the resulting group generator is: 
\begin{eqnarray}\label{final_group_generator}
\chi = (a_1+a_2r)\frac{\partial}{\partial r}+ (a_3+a_4t)\frac{\partial}{\partial t} + a_6\phi\frac{\partial}{\partial\phi}+(a_6-a_2)w\frac{\partial}{\partial w}.
\end{eqnarray}
In the above, the diffusion coefficient and its $r$-derivative are treated independently of each other for the purposes of invariance analysis. Therefore, the relationship between $D$ and $D_r$ must also be invariant under $\chi$. This is a straightforward calculation and it is shown in Appendix~A.  

Otherwise, Eqs.~(\ref{D-eq})-(\ref{a5-eqn}) include several quasi-linear first-order PDEs that will be solved using the Method of Characteristics. Equation~(\ref{final_group_generator}) may also be used to construct a system of similarity variables in order to start building a solution for Eq.~(\ref{timedependent_diff_eq}). However, the motivation of the current work is to understand the material properties and their connection to the presence or absence of symmetries. These results follow in Sec.~IV. 
\section{ANALYSIS OF RESULTS}

Equations~(\ref{t_r_balance})-(\ref{a5-eqn}) must be satisfied for Eq.~(\ref{eqn:mu1}) to be invariant under the groups of translation and scaling transformations generated by Eq.~(\ref{final_group_generator}). Practically, these determining equations amount to required values the constants $a_1 - a_8$ must assume, and associated functional forms assumed by the diffusion coefficient $D$ and macroscopic cross-section data encapsulated in $\Gamma$, respectively. Equations~(\ref{t_r_balance})-(\ref{a5-eqn}) contain several simple and intuitively obvious members:

\begin{itemize}
\item Equation~(\ref{a5-eqn}) indicates that $a_5 = 0$, or that translation in the scalar flux $\phi$ is never present. This result is obvious through inspection of Eq.~(\ref{eqn:diffeq_final}).

\item Equation~(\ref{t_r_balance}) indicates that $a_7 = 0$, or that translation in the spatial derivative of the scalar flux $w$ is also never present. This result is again obvious through inspection of Eq.~(\ref{eqn:w_def}). 

\item Equation~(\ref{phi_t_balance}) indicates that $a_8 = a_6 - a_2$, or that $w$ scales as $\phi / r$. Given the definition of $w$ appearing in Eq.~(\ref{eqn:w_def}), this result is also intuitively obvious.

\item Equation~(\ref{nD_Dequation}) contains three possibilites: $a_1 = 0$, $D_r =0$, or $a_1 = D_r =0$. The first conditions indicates scenarios where space translation symmetry is not present, while the second includes those where the functional form of the diffusion coefficient is only dependent on time. The last includes a scenario where space translation symmetry is not present with a diffusion coefficient with only time dependence.

\item Finally, Eq.~(\ref{nDa1equation}) contains three non-trivial possibilities: $n = 0$, $a_1 = 0$, or $n = a_1 = 0$. The first of these scenarios indicates that space translation symmetry may be present only in planar geometry, and the second indicates that space translation symmetry is not present in curvilinear geometries. The last includes scenarios where space translation symmetry is not present in planar geometry. 

\end{itemize}
The first three constraints above result in considerable simplfication of the remaining determining equations for $D$ and $\Gamma$. In particular, Eq. (\ref{D-eq}) becomes

\begin{equation}\label{reduced_diffusion_coeff_equation}
\left(a_1+a_2r\right)D_r + \left( a_3 +a_4t\right)D_t = (2a_2-a_4)D,
\end{equation}
and Eq.~(\ref{Dr-eq}) becomes
\begin{equation}\label{reduced_diffusion_coeff_deriv_equation}
\left(a_1+a_2r\right)D_{rr} + \left( a_3 +a_4t\right)D_{rt} = (a_2-a_4)D_r,
\end{equation}
which may also be obtained by taking the $r$-derivative of Eq. (\ref{reduced_diffusion_coeff_equation}). These two determining equations are thus mutually consistent, and only Eq.~(\ref{reduced_diffusion_coeff_equation}) will need to be solved to determine the functional form of $D$.

Equations~(\ref{Dr-eq}) and (\ref{gamma_eqn}) may be even further simplified and for each case, Eqs.~(\ref{nD_Dequation}) and (\ref{nDa1equation}) are satisfied according to the specific constraints under investigation. These constraints result in six unique cases which are summarized in Table \ref{SummaryTable}. 

\begin{itemize}
\item Case A: $n = 0$

With $n = 0$, Eqs. (\ref{gamma_eqn}) and (\ref{reduced_diffusion_coeff_equation}) are unchanged. Therefore, the solution for the diffusion coefficient is given by

\begin{equation}
D(r,t) = \left(a_3 +a_4t\right)^{\frac{2a_2}{a_4} -1}G\left[\left(r+\frac{a_1}{a_2}\right)\left(a_3+a_4t\right)^{-\frac{a_2}{a_4}}\right]
\end{equation}
where $G$ is an arbitrary function of the indicated argument.

In addition, Eq. (\ref{gamma_eqn}) may be solved using the Method of Characteristics to yield

\begin{equation}\label{gamma_solution}
\Gamma(r,t) = \left(a_3 +a_4t\right)^{-1}F\left[\left(r+\frac{a_1}{a_2}\right)\left(a_3+a_4t\right)^{-\frac{a_2}{a_4}}\right]
\end{equation}
where $F$ is another arbitrary function of the indicated argument. 
\item Case B: $a_1 = 0$

With $a_1 = 0$, Eqs. (\ref{gamma_eqn}) and (\ref{reduced_diffusion_coeff_equation}) become, respectively,
\begin{equation}\label{a1_zero_gamma}
a_2r\Gamma_r + (a_3+a_4t)\Gamma_t + a_4\Gamma =0,
\end{equation}
\begin{equation}\label{a1_zero_reduced_diffusion_coeff_equation}
a_2rD_r + \left( a_3 +a_4t\right)D_t = (2a_2-a_4)D.
\end{equation}

Again using the Method of Characteristics, the solutions to Eqs. (\ref{a1_zero_gamma}) and (\ref{a1_zero_reduced_diffusion_coeff_equation}) are, respectively,
\begin{equation}\label{a1_zero_gamma_solution}
\Gamma(r,t) = \left(a_3 +a_4t\right)^{-1}F\left[r\left(a_3+a_4t\right)^{-\frac{a_2}{a_4}}\right],
\end{equation}
and
\begin{equation}\label{a1_zero_diffusion_solution}
D(r,t) = \left(a_3 +a_4t\right)^{\frac{2a_2}{a_4} -1}G\left[r\left(a_3+a_4t\right)^{-\frac{a_2}{a_4}}\right],
\end{equation}
where $F$ and $G$ are again arbitrary functions of the indicated argument.

\item Case C: $n = a_1 = 0$

With $n = a_1 = 0$, Eqs. (\ref{gamma_eqn}) and (\ref{reduced_diffusion_coeff_equation}) become Eqs. (\ref{a1_zero_gamma}) and (\ref{a1_zero_reduced_diffusion_coeff_equation}), respectively.

This case is identical to the case $a_1 = 0$ considered above and will yield Eqs. (\ref{a1_zero_gamma_solution}) and (\ref{a1_zero_diffusion_solution}) as solutions.

\item Case D: $n = D_r = 0$

With $n=D_r = 0$, Eqs. (\ref{gamma_eqn}) and (\ref{reduced_diffusion_coeff_equation}) become, respectively,
\begin{equation}\label{n_Dr_zero_gamma}
\left(a_1 + a_2r\right)\Gamma_r + \left(a_3+a_4t\right)\Gamma_t + a_4\Gamma =0,
\end{equation}
\begin{equation}\label{n_Dr_zero_reduced_diffusion_coeff_equation}
\left( a_3 +a_4t\right)D_t = \left(2a_2-a_4\right)D.
\end{equation}

Again using the Method of Characteristics, the solution to Eq. (\ref{n_Dr_zero_gamma}) is (\ref{gamma_solution}), while the solution to Eq. (\ref{n_Dr_zero_reduced_diffusion_coeff_equation}) is
\begin{equation}\label{n_Dr_diffusion_solution}
D(t) = C (a_3+a_4t)^{\frac{2a_2}{a_4}-1},
\end{equation}
where $C$ is an arbitrary constant, and is dependent on initial conditions.
\item Case E: $D_r = a_1 =0$

With $D_r = a_1 = 0$, Eqs. (\ref{gamma_eqn}) and (\ref{reduced_diffusion_coeff_equation}) become Eqs. (\ref{a1_zero_gamma}) and (\ref{n_Dr_zero_reduced_diffusion_coeff_equation}), respectively. This case will yield Eqs. (\ref{a1_zero_gamma_solution}) and (\ref{n_Dr_diffusion_solution}) as solutions.

\item Case F: $a_1 = n = D_r =0$

With $a_1 = n = D_r = 0$, Eqs. (\ref{gamma_eqn}) and (\ref{reduced_diffusion_coeff_equation}) become Eqs. (\ref{a1_zero_gamma}) and (\ref{n_Dr_zero_reduced_diffusion_coeff_equation}), respectively. This case will yield Eqs. (\ref{a1_zero_gamma_solution}) and (\ref{n_Dr_diffusion_solution}) as solutions.
\end{itemize}
\begin{table}[h!]
\centering
\hspace*{-2.5cm}
\begin{tabular}{|c|c||p{7cm}|p{7cm}|}
\hline
\multicolumn{4}{|c|}{Material Property Functions} \\
\hline
\cellcolor[HTML]{000000}&Constraint(s)& Diffusion Constant $D$ & Gamma Coefficient $\Gamma$\\
\hline
\hline
Case A&$n=0$&\begin{eqnarray}
D(r,t) &=&\left(a_3 + a_4t \right)^{\frac{2a_2}{a_4}-1}\nonumber\\
&&\times G\left[\left(r+\frac{a_1}{a_2}\right)\left(a_3+a_4t\right)^{-\frac{a_2}{a_4}}\right]\nonumber
\end{eqnarray}& \begin{eqnarray}
\Gamma(r,t) &=& \left(a_3+a_4t \right)^{-1}\nonumber\\
&&\times F\left[\left(r + \frac{a_1}{a_2}  \right)\left(a_3+a_4t \right)^{\frac{a_2}{a_4}} \right]\nonumber
\end{eqnarray}  \\
\hline
Case B&$a_1 =0$& \begin{eqnarray}
D(r,t) &=& \left(a_3+a_4\right)^{\frac{2a_2}{a_4}-1}\nonumber\\
&&\times G\left[r\left(a_3+a_4t\right)^{-\frac{a_2}{a_4}}\right]\nonumber
\end{eqnarray}&\begin{eqnarray}
\Gamma(r,t) &=&\left(a_3+a_4t\right)^{-1}\nonumber\\
&&\times F\left[r\left(a_3+a_4t\right)^{-\frac{a_2}{a_4}}\right]\nonumber
\end{eqnarray} \\
\hline
Case C&$n=a_1 =0$& \begin{eqnarray}
D(r,t) &=& \left(a_3+a_4\right)^{\frac{2a_2}{a_4}-1}\nonumber\\
&&\times G\left[r\left(a_3+a_4t\right)^{-\frac{a_2}{a_4}}\right]\nonumber
\end{eqnarray}&\begin{eqnarray}
\Gamma(r,t) &=&\left(a_3+a_4t\right)^{-1}\nonumber\\
&&\times F\left[r\left(a_3+a_4t\right)^{-\frac{a_2}{a_4}}\right]\nonumber
\end{eqnarray} \\
\hline
Case D&$n=D_r=0$& \begin{eqnarray}
D(t) = C\left(a_3+a_4t\right)^{\frac{2a_2}{a_4}-1}\nonumber
\end{eqnarray} & \begin{eqnarray}
\Gamma(r,t) &=& \left(a_3+a_4t \right)^{-1}\nonumber\\
&&\times F\left[\left(r + \frac{a_1}{a_2}  \right)\left(a_3+a_4t \right)^{\frac{a_2}{a_4}} \right]\nonumber
\end{eqnarray}  \\
\hline
Case E&$D_r = a_1 =0$&\begin{eqnarray}
D(t) = C\left(a_3+a_4t\right)^{\frac{2a_2}{a_4}-1}\nonumber
\end{eqnarray} & \begin{eqnarray}
\Gamma(r,t) &=&\left(a_3+a_4t\right)^{-1}\nonumber\\
&&\times F\left[r\left(a_3+a_4t\right)^{-\frac{a_2}{a_4}}\right]\nonumber
\end{eqnarray} \\
\hline
Case F&$a_1 = n=D_r =0$&\begin{eqnarray}
D(t) = C\left(a_3+a_4t\right)^{\frac{2a_2}{a_4}-1}\nonumber
\end{eqnarray} & \begin{eqnarray}
\Gamma(r,t) &=&\left(a_3+a_4t\right)^{-1}\nonumber\\
&&\times F\left[r\left(a_3+a_4t\right)^{-\frac{a_2}{a_4}}\right]\nonumber
\end{eqnarray} \\
\hline
\end{tabular}
\caption{The summary of Cases A-F is incorperated in the above table. Note, that $F$ and $G$ are arbitrary functions and $C$ is an arbitrary constant that depends on initial conditions.}
\label{SummaryTable}
\end{table}
\section{DISCUSSION AND CONCLUSIONS}

The preceding study has demonstrated the application of symmetry analysis techniques in an EDS setting to the monoenergetic neutron diffusion equation with space- and time-variable material properties. The principal outcomes of this analysis are space-time functional forms the diffusion coefficient $D$ and cross-section data encapsulated in $\Gamma$ must satisfy to enable the existence of general translation and scaling symmetries (in addition to other, simpler constraints that may be intuitively obvious). Many of these formulae appearing in Sec. IV contain arbitrary functions of the space and time variables, thus enabling the existence of translation and scaling phenomena under a wide variety of scenarios (many of which are expected to be physically relevant). Some notable features along these lines include:
\begin{itemize}
\item Space translation symmetry is available only in planar geometry, and as a result the admissible functional forms for $D$ and $\Gamma$ are somewhat broadened as compared to cases in curvilinear geometry, or others where space translations are otherwise ignored (i.e., $a_1 = 0$). This result is consistent with previous studies of the 1D compressible flow equations expressed in terms of arbitrary 1D curvilinear coordinates; for these and other equations, moving to 1D planar geometry typically enhances the rank of the corresponding Lie algebra~\cite{axford_noh}. 
\item A class of space-time scaling and time-translation symmetries exists even in curvilinear geometry when $D_r = 0$. These symmetries are notable in that while the restriction on $D$ is somewhat severe, there are few analogous constraints on $\Gamma$. This outcome suggests the presence of scenarios where the $r$-dependence of $D$ and $\Gamma$ need not necessarily be closely linked. Given that both these functions depend principally on nuclear cross-section data, the implications of these symmetries represent a potentially interesting avenue of future study. 
\item The admissible functional forms of both $D$ and $\Gamma$ corresponding to various special cases can be constructed from the generalized results appearing in Sec. IV. For example, further restrictions on $D$ and $\Gamma$ may be derived for cases where the constants $a_i$ assume specific values (including zero, in many cases, indicating the absence of the indicated symmetry). 
\item Many commonly encountered special cases can also be made to conform to the required functional forms of $D$ and $\Gamma$; for example, requiring that these functions be either constant or some other simple function of space and time in turn sets required values for the constants $a_i$. If at least one of these constants is then found to be non-zero, the results of this study directly produce the correspondence between a problem's physical conditions and the subset of symmetries they enable. 
\end{itemize}
Indeed, the outcomes of this work are closely related to those of Ramsey et al.~\cite{ramsey_tellez}, who have noted that the usual space-time separable solutions of the monoenergetic neutron diffusion equation have well-defined ties to translational and scaling symmetry groups. These results are implicitly encoded in the current study, which presumably includes broader classes of space-time separable and other solutions with special physical properties.

Along these lines, future studies stemming from this work may include:

\begin{itemize}
\item Constructing solutions of Eq. (\ref{eqn:diffeq_final}) using the translation and scale-invariant properties discovered in this work, and explicitly connecting them back to those of Ramsey et al. or others,

\item Determining the general symmetry groups admitted by Eq. (\ref{eqn:diffeq_final}), and their associated solutions,

\item Performing similar analyses for multi-energy group, multi-region models, $P_n$ or $S_n$ equations of various orders, or the integro-differential neutron transport equation in a variety of corrdinate sustems (e.g., generalized curvilinear coordinates).
\end{itemize}

As noted in Sec. I, the rigorous application of symmetry analysis techniques to the neutron transport equation or its surrogates appears to be largely unexplored to date, and thus represents fertile ground for the development of various new solutions and the development of physical intuition. 

\section{ACKNOWLEDGEMENT}
This work was performed under the auspices of the United
States Department of Energy by Los Alamos National Security, LLC, at Los Alamos National Laboratory under contract
DE-AC52-06NA25396.

The authors would like to thank E. J. Albright, R. Baty, Z. Boyd, P. Jaegers, E. Schmidt, J. Schmidt, and J. Tellez for their valuable insights. JFG would like to thank Prof. R. F. Lebed (Arizona State University) for proofreading and overseeing this work. 

\appendix
\numberwithin{equation}{section}
\section*{Appendix}
\section{Diffusion Completeness}

In Sec. III the diffusion coefficient $D$ and its $r$-derivative are treated as independent of each other for the purposes of invariance analysis. For completeness, it must also be true that the definition of $D_r$ as an $r$-derivative is invariant under the same group of transformations generated by $\chi$ appearing in Eq. (\ref{group_generator}). To begin, a two-form corresponding to the definition of $D_r$, namely 
\begin{equation}\label{Drpartial}
D_r = \frac{\partial D}{\partial r}
\end{equation}
is given by 
\begin{equation}\label{eqn:mu_3}
\mu_3 = D_r~dr\wedge dt - dD\wedge dt.
\end{equation}
where it can be shown that Eqs. (\ref{Drpartial}) and (\ref{eqn:mu_3}) are equivalent by sectioning and annulling as done in Section II. B. 

The condition that Eq. (\ref{eqn:mu_3}) is invariant is given by 
\begin{eqnarray}\label{chi_mu3}
\chi\mu_3 = \big[\left(a_1 + a_2 r\right)D_{rr}+\left(a_3+a_4t\right)D_{rt}\big]dr\wedge dt 
+ D_r\left(a_2+a_4\right)dr\wedge dt-a_4 dD\wedge dt\nonumber\\
-d\big[\left(a_1 + a_2 r\right)D_{r}+\left(a_3+a_4t\right)D_{t}\big]
\end{eqnarray}
where the condition $\mu_3=0$ will be implemented later. To further simplify this condition, the identity
\begin{equation}\label{transitive_identity}
d\frac{\partial D}{\partial q} = \frac{\partial}{\partial q}\bigg[dD\bigg],
\end{equation} 
may be applied, where $q$ is any coordinate. Moreover,
\begin{equation}\label{Dr_transitive}
\frac{\partial}{\partial r}dD = \frac{\partial^2 D}{\partial r^2}dr + \frac{\partial D}{\partial r}d\bigg[\frac{\partial r}{\partial r}\bigg] + \frac{\partial^2 D}{\partial r\partial t}dt + \frac{\partial D}{\partial t}d\bigg[\frac{\partial t}{\partial r}\bigg],
\end{equation}
\begin{equation}\label{Dt_transitive}
\frac{\partial}{\partial t}dD = \frac{\partial^2 D}{\partial r\partial t}dr + \frac{\partial D}{\partial r}d\bigg[\frac{\partial t}{\partial r}\bigg] + \frac{\partial^2 D}{\partial t^2}dt + \frac{\partial D}{\partial t}d\bigg[\frac{\partial t}{\partial t}\bigg].
\end{equation}
Since $t, r$ are independent of each other and the derivative of a constant is zero, Eqs. (\ref{Dr_transitive}) and (\ref{Dt_transitive}) become, respectively,
\begin{eqnarray}\label{Dr_transitive_reduced}
\frac{\partial}{\partial r}dD = \frac{\partial^2 D}{\partial r^2}dr + \frac{\partial^2 D}{\partial r\partial t}dt,
\end{eqnarray}
\begin{eqnarray}\label{Dt_transitive_reduced}
\frac{\partial}{\partial t}dD = \frac{\partial^2 D}{\partial r\partial t}dr + \frac{\partial^2 D}{\partial t^2}dt.
\end{eqnarray}

Finally, with Eqs. (\ref{transitive_identity}), (\ref{Dr_transitive}), (\ref{Dt_transitive}) and $\mu_3 =0$, Eq. (\ref{chi_mu3}) becomes
\begin{eqnarray}\label{chi_mu3_solved}
\chi\mu_3 &=& \big[\left(a_1 + a_2 r\right)D_{rr}+\left(a_3+a_4t\right)D_{rt}\big]dr\wedge dt + D_r\left(a_2+a_4\right)dr\wedge dt-a_4  D_r dr\wedge dt\nonumber\\
&&-\bigg\{ D_r a_2dr + (a_1 +a_2r)\big[D_{rr}dr + D_{rt}dt\big]
+ a_4D_tdt+(a_3+a_4t)\big[D_{rt}dr +D_{tt}dt\big] \wedge dt\bigg\}\nonumber\\
&=&0\nonumber\\
\rightarrow 
\chi\mu_3 &=& \big[\left(a_1 + a_2 r\right)D_{rr}+\left(a_3+a_4t\right)D_{rt}\big]dr\wedge dt
+ D_r\left(a_2+a_4\right)dr\wedge dt-D_r(a_2+a_4) dr\wedge dt\nonumber\\
&&-\big[(a_1 +a_2r)D_{rr}+(a_3+a_4t)D_{rt}\big] dr\wedge dt=0,
\end{eqnarray}
which is identically satisfied as indicated. Therefore, the closure relation given by Eq. (\ref{Drpartial}) is also invariant under the group transformations generated by $\chi$, as expected.

\end{document}